\documentclass[11pt]{article}
\usepackage{amsfonts}
\usepackage{amsmath}
\usepackage{amssymb}
\usepackage{graphicx}
\usepackage[dvipsnames,usenames]{color}
\usepackage[pdfstartview=FitH ]{hyperref}
\usepackage[hyperpageref]{backref}
\begin{document}

\baselineskip 18pt
\def\o{\over}
\def\e{\varepsilon}
\title{\Large\bf  Kloosterman\ \ Sums\ \ with\ \ Multiplicative\ \
Coefficients}
\author{Ke\ Gong\ and\ Chaohua\ Jia}
\date{}
\maketitle {\small \noindent {\bf Abstract.} Let $f(n)$ be a
multiplicative function satisfying $|f(n)|\leq 1$, $q$ $(\leq N^2)$
be a positive integer and $a$ be an integer with $(a,\,q)=1$. In
this paper, we shall prove that
$$
\sum_{\substack{n\leq N\\ (n,\,q)=1}}f(n)e({a\bar{n}\o q})\ll
\sqrt{\tau(q)\o q}N\log\log(6N)+q^{{1\o 4}+{\e\o 2}}N^{1\o
2}(\log(6N))^{1\o 2}+{N\o \sqrt{\log\log(6N)}},
$$
where $\bar{n}$ is the multiplicative inverse of $n$ such that
$\bar{n}n\equiv 1\,({\rm mod}\,q),\,e(x)=\exp(2\pi ix),\,\tau(q)$ is
the divisor function. }

\vskip.3in
\noindent{\bf 1. Introduction}

Let $\mu(n)$ be the M\"obius function, $q$ be a positive integer and
$a$ be an integer with $(a,\,q)=1$. In 1988, D. Hajela, A.
Pollington and B. Smith [8] proved that
$$
\sum_{\substack{n\leq N\\ (n,\,q)=1}}\mu(n)e({a\bar{n}\o q})\ll_\e
Nq^\e\Bigl({(\log N)^{5\o 2}\o q^{1\o 2}}+{q^{3\o 10}(\log N)^{11\o
5}\o N^{1\o 5}}\Bigr),
$$
where $\bar{n}$ is the multiplicative inverse of $n$ such that
$\bar{n}n\equiv 1\,({\rm mod}\,q)$, $e(x)=\exp(2\pi ix)$ and $\e$ is
a sufficiently small positive constant. This estimate is nontrivial
for $(\log N)^{5+10\e}\ll q\ll N^{{2\o 3}-3\e}$.

Later, P. Deng [4], G. Wang and Z. Zheng [9] independently improved
the above estimate to
$$
\sum_{\substack{n\leq N\\ (n,\,q)=1}}\mu(n)e({a\bar{n}\o q})\ll
N\tau(q)\Bigl({(\log N)^{5\o 2}\o q^{1\o 2}}+{q^{1\o 5}(\log
N)^{13\o 5}\o N^{1\o 5}}\Bigr),
$$
where $\tau(q)$ is the divisor function, which is nontrivial for
$(\log N)^{5+\e}\ll q\ll N^{1-\e}$. It was stated in [4] that under
the Generalized Riemann Hypothesis, one can get
$$
\sum_{\substack{n\leq N\\ (n,\,q)=1}}\mu(n)e({a\bar{n}\o q})\ll_\e
q^{1\o 2}N^{{1\o 2}+\e}.
$$

We also mention some progress on the relative topic. In 1998, E.
Fouvry and P. Michel [6] proved that if $q$ is a prime number,
$g(x)={P(x)\o Q(x)}$ is any rational function with $P(x)$ and $Q(x)$
relatively prime monic polynomials in ${\Bbb Z}[x]$, then for $1\leq
N\leq q$, one has
$$
\sum_{\substack{p\leq N\\ (Q(p),\,q)=1}}e({g(p)\o q})\ll_\e q^{{3\o
16}+\e}N^{25\o 32},
$$
where $p$ runs through prime numbers, the implied constant also
depends on the degrees of $P$ and $Q$. This estimate is nontrivial
for $N\leq q\ll N^{{7\o 6}-7\e}$. It was stated in [6] that the same
method can produce
$$
\sum_{\substack{n\leq N\\ (Q(n),\,q)=1}}\mu(n)e({g(n)\o q})\ll_\e
q^{{3\o 16}+\e}N^{25\o 32}
$$
for the prime number $q$ and $1\leq N\leq q$. Some further results
can be found in [5].

In 2011, E. Fouvry and I. E. Shparlinski [7] proved that for
$(a,\,q)=1$ and $N^{3\o 4}\leq q\leq N^{4\o 3}$, one has
$$
\sum_{\substack{N<p\leq 2N\\ (p,\,q)=1}}e({a\bar{p}\o q})\ll_\e
q^\e(q^{1\o 4}N^{2\o 3}+N^{15\o 16}),
$$
which is nontrivial for $N^{3\o 4}\leq q\ll N^{{4\o 3}-6\e}$. They
also proved that if $(a,\,q)=1$, then
$$
\sum_{\substack{N<p\leq 2N\\ (p,\,q)=1}}e({a\bar{p}\o q})\ll
N\Bigl(\tau^{1\o 2}(q){(\log N)^2\o q^{1\o 2}}+\tau(q){q^{1\o
4}(\log N)^{3\o 2}\o N^{1\o 5}}\Bigr),
$$
which is nontrivial for $(\log N)^{6+\e}\ll q\ll N^{{4\o 5}-\e}$. In
2012, R. C. Baker [1] gave improvement under some conditions.

When the first author visited the University of Montreal, Professor
A. Granville suggested him to study the general sum
$$
\sum_{\substack{n\leq N\\ (n,\,q)=1}}f(n)e({a\bar{n}\o q}), \eqno
(1.1)
$$
where $f(n)$ is a multiplicative function satisfying $|f(n)|\leq 1$.

In this paper, we shall apply the method in Section 2 of [3], which
is called as the finite version of Vinogradov's inequality, to give
a nontrivial estimate for the sum in (1.1) when $q$ is in a suitable
range.

{\bf Theorem}. Let $f(n)$ be a multiplicative function satisfying
$|f(n)|\leq 1$, $q$ $(\leq N^2)$ be a positive integer and $a$ be an
integer with $(a,\,q)=1$. Then we have
\begin{align*}
\qquad\qquad &\sum_{\substack{n\leq N\\
(n,\,q)=1}}f(n)e({a\bar{n}\o q})\ll
\sqrt{\tau(q)\o q}N\log\log(6N)\qquad\qquad\qquad\quad\ \ \ (1.2)\\
&\qquad\qquad +q^{{1\o 4}+{\e\o 2}}N^{1\o 2}(\log(6N))^{1\o 2}+{N\o
\sqrt{\log\log(6N)}}.
\end{align*}

The estimate in (1.2) is nontrivial for
$$
(\log\log(6N))^{2+\e}\ll q\ll N^{2-5\e}.
$$
In a private communication, Ping Xi remarked that when $q$ is a
prime number, if Lemma 2 below is replaced by Theorem 16 in [2],
then the upper bound in the above nontrivial range can be extended
to $q\ll N^A$, where $A$ is any given large constant.

Throughout this paper, we assume that $N$ is sufficiently large and
set
\begin{align*}
\qquad\qquad\quad &d_0=\sqrt{\log\log(6N)},\quad
D_0=e^{d_0}=\exp(\sqrt{\log\log(6N)}),
\qquad (1.3)\\
&d_1=d_0^2=\log\log(6N),\quad D_1=e^{d_1}=\log(6N).
\end{align*}
Let $p$ denote a prime number, $\tau(q)$ denote the divisor
function, $\e$ be a sufficiently small positive constant.

\vskip.3in
\noindent{\bf 2. Some preliminaries}

Write
\begin{align*}
\qquad\quad\ S&=\{n:\,1\leq n\leq N,\, n{\rm\ has\ a\ prime\ factor\
in\ }[D_0,\,
D_1)\},\qquad\quad (2.1)\\
T&=\{n:\,1\leq n\leq N,\, n{\rm\ has\ no\ prime\ factor\ in\
}[D_0,\,D_1)\}.
\end{align*}

{\bf Lemma 1}. We have
$$
|T|\ll{N\o \sqrt{\log\log(6N)}}.
$$

Proof. Let
$$
P(N)=\prod_{D_0\leq p<D_1}p.
$$
We have
\begin{align*}
|T|&=\sum_{\substack{n\leq N\\ (n,\,P(N))=1}}1\\
&=\sum_{n\leq N}\sum_{d|(n,\,P(N))}\mu(d)\\
&=\sum_{d|P(N)}\mu(d)\sum_{\substack{n\leq N\\ d|n}}1\\
&=\sum_{d|P(N)}\mu(d)\Bigl({N\o d}+O(1)\Bigr)\\
&=N\sum_{d|P(N)}{\mu(d)\o d}+O\Bigl(2^{\pi(D_1)}\Bigr)\\
&=N\prod_{D_0\leq p<D_1}\Bigl(1-{1\o p}\Bigr)+O\Bigl(2^{2D_1\o \log
D_1}\Bigr)\\
&\ll N{\log D_0\o \log D_1}+O\Bigl(2^{2\log(6N)\o \log\log(6N)}\Bigr)\\
&\ll{N\o \sqrt{\log\log(6N)}}.
\end{align*}
Hence, Lemma 1 holds true.

By Lemma 1, we have
$$
\sum_{\substack{n\leq N\\ (n,\,q)=1}}f(n)e({a\bar{n}\o q})=
\sum_{\substack{n\leq N\\ n\in S\\ (n,\,q)=1}}f(n)e({a\bar{n}\o q})
+O\Bigl({N\o \sqrt{\log\log(6N)}}\Bigr). \eqno(2.2)
$$

Let
$$
P_r=\{p:\,e^r\leq p< e^{r+1}\},\qquad\quad{\rm if}\ \ [d_0]\leq
r\leq [d_1]. \eqno (2.3)
$$
Then
$$
\bigcup_{r=[d_0]+1}^{[d_1]-1}P_r\subseteq\{p:\,D_0\leq
p<D_1\}\subseteq\bigcup_{r=[d_0]}^{[d_1]}P_r.
$$
The prime number theorem yields
$$
|P_r|\ll{e^r\o r}. \eqno (2.4)
$$

Write
\begin{align*}
S'&=\{n:\,1\leq n\leq N,\, n{\rm\ has\ a\ prime\ factor\ in\
}\bigcup_{r=[d_0]}^{[d_1]}P_r\},\\
S''&=\{n:\,1\leq n\leq N,\, n{\rm\ has\ a\ prime\ factor\ in\
}\bigcup_{r=[d_0]+1}^{[d_1]-1}P_r\}.
\end{align*}
Then
$$
S''\subseteq S\subseteq S'.
$$
Hence,
\begin{align*}
|S\backslash S''|\leq|S'\backslash S''|&\ll\sum_{p\in P_{[d_0]}}{N\o
p}+\sum_{p\in P_{[d_1]}}{N\o p}\\
&\ll N\Bigl({|P_{[d_0]}|\o e^{[d_0]}}+{|P_{[d_1]}|\o e^{[d_1]}}\Bigr)\\
&\ll{N\o d_0}={N\o \sqrt{\log\log(6N)}}.
\end{align*}

We note that
\begin{align*}
&\ |\{n:\,1\leq n\leq N,\, n{\rm\ has\ at\ least\ two\ prime\
factors\ in\ the}\\
&\qquad\qquad\quad{\rm same\ one\ of\ }P_r's\,([d_0]+1\leq r\leq
[d_1]-1)\}|\\
&\ll\sum_{r=[d_0]+1}^{[d_1]-1}\sum_{p\in P_r}\sum_{p'\in P_r}{N\o
pp'}\\
&\ll N\sum_{r=[d_0]+1}^{[d_1]-1}\Bigl({|P_r|\o e^r}\Bigr)^2\\
&\ll N\sum_{r=[d_0]+1}^{[d_1]-1}{1\o r^2}\\
&\ll{N\o d_0}={N\o \sqrt{\log\log(6N)}}.
\end{align*}
Therefore for
\begin{align*}
&S'''=\{n:\,1\leq n\leq N,\, n{\rm\ has\ exact\ one\ prime\ factor}\\
&\qquad\qquad\qquad{\rm in\ one\ of\ }P_r's\,([d_0]+1\leq r\leq
[d_1]-1)\},
\end{align*}
we have
$$
S'''\subseteq S''
$$
and
$$
|S''\backslash S'''|\ll{N\o \sqrt{\log\log(6N)}}.
$$

The set $S'''$ can be decomposed as
$$
S'''=\bigcup_{r=[d_0]+1}^{[d_1]-1}S_r, \eqno (2.5)
$$
where
\begin{align*}
&\qquad\  S_r=\{n:\,1\leq n\leq N,\, n{\rm\ has\ exact\ one\
prime\ factor\ in\ }P_r\qquad\quad (2.6)\\
&\qquad\qquad\qquad\qquad {\rm and\ has\ no\ prime\ factor\ in\
}\bigcup_{i<r}P_i\}.
\end{align*}
By the prime number theorem, it is easy to see that each
$S_r(r=[d_0]+1,\,\cdots,\,[d_1]-1)$ is not empty. The sets $S_r$ are
disjoint from each other. Every element $n\in S_r$ can be written in
exact one way as
$$
n=py,\eqno (2.7)
$$
where $p\in P_r,\,y$ has no prime factor in $\bigcup_{i\leq
r}P_i,\,py\leq N$.

From the above discussion, we get
\begin{align*}
\qquad &\ \,\sum_{\substack{n\leq N\\ (n,\,q)=1}}f(n)e({a\bar{n}\o q})\\
&=\sum_{\substack{n\leq N\\ n\in S'''\\ (n,\,q)=1}}f(n)e({a\bar{n}\o
q}) +O\Bigl({N\o \sqrt{\log\log(6N)}}\Bigr)\\
&=\sum_{r=[d_0]+1}^{[d_1]-1}\sum_{\substack{n\leq N\\ n\in S_r\\
(n,\,q)=1}}f(n)e({a\bar{n}\o q}) +O\Bigl({N\o \sqrt{\log\log(6N)}}\Bigr)\\
&=\sum_{r=[d_0]+1}^{[d_1]-1}\sum_{\substack{e^r\leq p< e^{r+1}\\
(p,\,q)=1}}\sum_{\substack{y\leq {N\o p}\\ y{\rm\tiny\ has\ no\
prime\ factor\ in\  }\bigcup_{i\leq r}P_i\\ (y,\,q)=1}}f(py)
e({a\bar{p}\bar{y}\o q})\qquad\ \ \ (2.8)\\
&\qquad\qquad\qquad +O\Bigl({N\o \sqrt{\log\log(6N)}}\Bigr)
\end{align*}
\begin{align*}
&=\sum_{r=[d_0]+1}^{[d_1]-1}\sum_{\substack{y\leq {N\o e^r}\\
y{\rm\tiny\ has\ no\ prime\ factor\ in\ }\bigcup_{i\leq r}P_i\\
(y,\,q)=1}}f(y)\sum_{\substack{e^r\leq p< e^{r+1}\\
p\leq{N\o y}\\ (p,\,q)=1}}f(p)e({a\bar{p}\bar{y}\o q})\\
&\qquad\qquad\qquad+O\Bigl({N\o \sqrt{\log\log(6N)}}\Bigr)\\
&\ll\sum_{r=[d_0]+1}^{[d_1]-1}\sum_{\substack{y\leq {N\o e^r}\\
(y,\,q)=1}}\Bigl|\sum_{\substack{e^r\leq p< e^{r+1}\\
p\leq{N\o y}\\ (p,\,q)=1}}f(p)e({a\bar{p}\bar{y}\o q})\Bigr| +{N\o
\sqrt{\log\log(6N)}}.
\end{align*}

Let
$$
Y={N\o e^r}. \eqno (2.9)
$$
We shall estimate the sum
$$
{\sum}_1=\sum_{\substack{y\leq Y\\ (y,\,q)=1}}\Bigl|\sum_{\substack
{e^r \leq p< e^{r+1}\\ p\leq{N\o y}\\ (p,\,q)=1}}f(p)e({a\bar{p}
\bar{y}\o q})\Bigr|. \eqno (2.10)
$$

{\bf Lemma 2}. For the positive integer $q$ and the integer $b$, we
have
$$
\sum_{\substack{X<n\leq Z\\ (n,\,q)=1}}e({b\bar{n}\o q})\ll
\Bigl({Z-X\o q}+1\Bigr)(b,\,q)+q^{{1\o 2}+\e}. \eqno (2.11)
$$

Proof. Lemma 2.1 in [7] states that
\begin{align*}
\sum_{\substack{X<n\leq Z\\ (n,\,q)=1}}e({b\bar{n}\o
q})&\ll\mu^2\Bigl({q\o (b,\,q)}\Bigr)\Bigl({Z-X\o
q}+1\Bigr)\cdot{\varphi(q)\o \varphi\Bigl({q\o (b,\,q)}\Bigr)}\\
&\ \ +\tau(q)\tau((b,\,q))\log(2q)q^{1\o 2}.
\end{align*}
Then the bounds
\begin{align*}
{\varphi(q)\o \varphi({q\o (b,\,q)})}&=q\prod_{p|q}\Bigl(1-{1\o
p}\Bigr)\cdot\Bigl({q\o (b,\,q)}\prod_{p|{q\o (b,\,q)}}\Bigl(1-{1\o
p}\Bigr)\Bigr)^{-1}\\
&=(b,\,q)\prod_{\substack{p|q\\ p\not\,|{q\o (b,\,q)}}}\Bigl(1-{1\o
p}\Bigr)\leq (b,\,q)
\end{align*}
and
$$
\tau(q)\ll q^{\e\o 4}
$$
produce the conclusion in Lemma 2.

\vskip.3in
\noindent{\bf 3. The proof of Theorem}

By Cauchy's inequality,
$$
{\sum}_1\leq Y^{1\o 2}\Bigl(\sum_{\substack{y\leq Y\\
(y,\,q)=1}}\Bigl|\sum_{\substack {e^r \leq p< e^{r+1}\\ p\leq{N\o
y}\\ (p,\,q)=1}}f(p)e({a\bar{p} \bar{y}\o q})\Bigr|^2\Bigr)^{1\o 2}.
\eqno (3.1)
$$
An application of Lemma 2 to
$$
{\sum}_2=\sum_{\substack{y\leq Y\\ (y,\,q)=1}}\Bigl|\sum_{\substack
{e^r \leq p< e^{r+1}\\ p\leq{N\o y}\\ (p,\,q)=1}}f(p)e({a\bar{p}
\bar{y}\o q})\Bigr|^2
$$
produces
\begin{align*}
\ \ \ {\sum}_2&=\sum_{\substack{y\leq Y\\ (y,\,q)=1}}\sum_{\substack
{e^r \leq p_1< e^{r+1}\\ p_1\leq{N\o y}\\
(p_1,\,q)=1}}\sum_{\substack
{e^r \leq p_2< e^{r+1}\\ p_2\leq{N\o y}\\
(p_2,\,q)=1}}f(p_1)\overline{f(p_2)}e({a(\bar{p}_1-\bar{p}_2)
\bar{y}\o q})\\
&=\sum_{\substack {e^r \leq p_1<
e^{r+1}\\(p_1,\,q)=1}}\sum_{\substack {e^r \leq p_2< e^{r+1}\\
(p_2,\,q)=1}}f(p_1)\overline{f(p_2)}\sum_{\substack{y\leq Y\\
y\leq{N\o \max(p_1,\,p_2)}\\(y,\,q)=1}}e({a(\bar{p}_1-\bar{p}_2)
\bar{y}\o q})\\
&\ll\sum_{\substack {e^r \leq p_1<
e^{r+1}\\(p_1,\,q)=1}} \sum_{\substack {e^r \leq p_2< e^{r+1}\\
(p_2,\,q)=1}}\Bigl|
\sum_{\substack{y\leq Y\\
y\leq{N\o \max(p_1,\,p_2)}\\(y,\,q)=1}}e({a(\bar{p}_1-\bar{p}_2)
\bar{y}\o q})\Bigr|\qquad\qquad  (3.2)\\
&\ll\sum_{\substack{e^r\leq p<e^{r+1}\\ (p,\,q)=1}}Y+\sum_{\substack
{e^r \leq p_1<
e^{r+1}\\(p_1,\,q)=1}} \sum_{\substack {e^r \leq p_2< e^{r+1}\\
(p_2,\,q)=1\\ p_2\ne p_1}}\Bigl|
\sum_{\substack{y\leq Y\\
y\leq{N\o \max(p_1,\,p_2)}\\(y,\,q)=1}}e({a(\bar{p}_1-\bar{p}_2)
\bar{y}\o q})\Bigr|\\
&\ll Ye^r+\sum_{\substack {e^r \leq p_1< e^{r+1}\\(p_1,\,q)=1}}
\sum_{\substack {e^r \leq p_2< e^{r+1}\\
(p_2,\,q)=1\\ p_2\ne p_1}}\Bigl(\Bigl({Y\o
q}+1\Bigr)(a(\bar{p}_1-\bar{p}_2),\,q)+q^{{1\o 2}+\e}\Bigr)\\
&\ll Ye^r+\Bigl({Y\o q}+1\Bigr)\sum_{\substack {e^r \leq p_1<
e^{r+1}\\(p_1,\,q)=1}}\sum_{\substack {e^r \leq p_2< e^{r+1}\\
(p_2,\,q)=1\\ p_2\ne p_1}}(\bar{p}_1-\bar{p}_2,\,q)+q^{{1\o
2}+\e}e^{2r}.
\end{align*}

We have
\begin{align*}
\qquad\ \ \ &\ \,\sum_{\substack {e^r \leq p_1<
e^{r+1}\\(p_1,\,q)=1}}
\sum_{\substack {e^r \leq p_2< e^{r+1}\\
(p_2,\,q)=1\\ p_2\ne p_1}}(\bar{p}_1-\bar{p}_2,\,q)\\
&=\sum_{k|q}k\sum_{\substack {e^r \leq p_1<
e^{r+1}\\(p_1,\,q)=1}}\sum_{\substack {e^r \leq p_2< e^{r+1}\\
(p_2,\,q)=1\\ p_2\ne p_1\\ (\bar{p}_1-\bar{p}_2,\,q)=k}} 1\qquad\qquad
\qquad\qquad\qquad\ (3.3)\\
&\leq \sum_{k|q}k\sum_{\substack {e^r \leq p_1<
e^{r+1}\\(p_1,\,q)=1}}\sum_{\substack {e^r \leq p_2< e^{r+1}\\
(p_2,\,q)=1\\ p_2\ne p_1\\ \bar{p}_2\equiv \bar{p}_1\,({\rm
mod}\,k)}}1\\
&=\sum_{k|q}k\sum_{\substack {e^r \leq p_1<
e^{r+1}\\(p_1,\,q)=1}}\sum_{\substack {e^r \leq p_2< e^{r+1}\\
(p_2,\,q)=1\\ p_2\ne p_1\\ p_2\equiv p_1\,({\rm mod}\,k)}}1.
\end{align*}
In the above sum, if $k\geq e^{r+1}$, then $p_2\equiv p_1\,({\rm
mod}\,k)$ and $p_1,\,p_2<e^{r+1}\Longrightarrow p_2=p_1$, which
contradicts the fact $p_2\ne p_1$. Hence, it follows that
\begin{align*}
&\ \,\sum_{k|q}k\sum_{\substack {e^r \leq p_1<
e^{r+1}\\(p_1,\,q)=1}}\sum_{\substack {e^r \leq p_2< e^{r+1}\\
(p_2,\,q)=1\\ p_2\ne p_1\\ p_2\equiv p_1\,({\rm mod}\,k)}}1\\
&=\sum_{\substack{k|q\\ k< e^{r+1}}}k\sum_{\substack {e^r \leq p_1<
e^{r+1}\\(p_1,\,q)=1}}\sum_{\substack {e^r \leq p_2< e^{r+1}\\
(p_2,\,q)=1\\ p_2\ne p_1\\ p_2\equiv p_1\,({\rm mod}\,k)}}1\\
&\leq\sum_{\substack{k|q\\ k<
e^{r+1}}}k\sum_{n_1<e^{r+1}}\sum_{\substack{n_2<e^{r+1}\\ n_2\equiv
n_1\,({\rm mod}\,k)}}1\\
&\ll\sum_{\substack{k|q\\ k< e^{r+1}}}k\cdot e^{r+1}\cdot{e^{r+1}\o
k}\\
&\ll\tau(q)e^{2r}.
\end{align*}
Thus we get the estimate
$$
\sum_{\substack {e^r \leq p_1<
e^{r+1}\\(p_1,\,q)=1}}\sum_{\substack {e^r \leq p_2< e^{r+1}\\
(p_2,\,q)=1\\ p_2\ne p_1}}(\bar{p}_1-\bar{p}_2,\,q)\ll\tau(q)e^{2r}.
\eqno (3.4)
$$

By the above discussion, we have
\begin{align*}
{\sum}_2&\ll Ye^r+\Bigl({Y\o q}+1\Bigr)\tau(q)e^{2r}
+q^{{1\o 2}+\e}e^{2r}\\
&\ll{\tau(q)\o q}Ye^{2r}+Ye^r+q^{{1\o 2}+\e}e^{2r}.
\end{align*}
It follows that
\begin{align*}
{\sum}_1&\ll Y^{1\o 2}\Bigl({\tau(q)\o
q}Ye^{2r}+Ye^r+q^{{1\o 2}+\e}e^{2r}\Bigr)^{1\o 2}\\
&\ll\sqrt{\tau(q)\o q}Ye^r+Ye^{r\o 2}+Y^{1\o 2}q^{{1\o
4}+{\e\o 2}}e^r\\
&\ll\sqrt{\tau(q)\o q}N+{N\o e^{r\o 2}}+q^{{1\o 4}+{\e\o 2}}N^{1\o
2}e^{r\o 2}.
\end{align*}
Applying this estimate to (2.8), we get
\begin{align*}
&\ \,\sum_{\substack{n\leq N\\ (n,\,q)=1}}f(n)e({a\bar{n}\o q})\\
&\ll\sum_{r=[d_0]+1}^{[d_1]-1}\Bigl(\sqrt{\tau(q)\o q}N+{N\o e^{r\o
2}}+q^{{1\o 4}+{\e\o 2}}N^{1\o
2}e^{r\o 2}\Bigr)+{N\o \sqrt{\log\log(6N)}}\\
&\ll\sqrt{\tau(q)\o q}N\log\log(6N)+{N\o
\exp({1\o 2}\sqrt{\log\log(6N)})}\\
&\qquad\qquad+q^{{1\o 4}+{\e\o 2}}N^{1\o 2}(\log(6N))^{1\o
2}+{N\o \sqrt{\log\log(6N)}}\\
&\ll\sqrt{\tau(q)\o q}N\log\log(6N)+q^{{1\o 4}+{\e\o 2}}N^{1\o
2}(\log(6N))^{1\o 2}+{N\o \sqrt{\log\log(6N)}}.
\end{align*}

So far the proof of Theorem is complete.

\vskip.3in
\noindent{\bf Acknowledgements}

The first author would like to thank Professor A. Granville for
giving him valuable suggestions. Both authors would like to thank
Professor P. Sarnak for his fascinating talks addressed in January
2014 in Jinan, which attracted their attention to the method in [3].
They would like to thank Professor Jianya Liu for kindly inviting
them to attend the conferment ceremony of honorary doctorate on
Professor P. Sarnak in Shandong University. They also would like to
thank Professor E. Fouvry and Ping Xi for their nice comments and
suggestions.

The first author is supported by the National Natural Science
Foundation of China (Grant No. 11126150) and the Natural Science
Foundation of the Education Department of Henan Province (Grant No.
2011A110003). The second author is supported by the National Natural
Science Foundation of China (Grant No. 11371344) and the National
Key Basic Research Program of China (Project No. 2013CB834202).


\bigskip

\

\

Ke Gong

Department of Mathematics, Henan University, Kaifeng, Henan 475004,
P. R. China

E-mail: kg@henu.edu.cn

\

Chaohua Jia

Institute of Mathematics, Academia Sinica, Beijing 100190, P. R.
China

Hua Loo-Keng Key Laboratory of Mathematics, Chinese Academy of
Sciences, Beijing 100190, P. R. China

E-mail: jiach@math.ac.cn


\vskip.6in

\begin{thebibliography}{9}

\bibitem{1} R. C. Baker, {\it Kloosterman sums with prime variable},
Acta Arith., {\bf 156}(2012), no.4, 351-372.

\bibitem{2} J. Bourgain and M. Z. Garaev, {\it Sumsets of
reciprocals in prime fields and multilinear Kloosterman sums},
available at http://arxiv.org/abs/1211.4184v1.

\bibitem{3} J. Bourgain, P. Sarnak and T. Ziegler, {\it Disjointness
of Moebius from horocycle flows}, From Fourier Analysis and Number
Theory to Radon Transforms and Geometry, 67-83, Developments in
Mathematics 28, Springer, New York, 2013.

\bibitem{4} P. Deng, {\it On Kloosterman sums with oscillating coefficients},
Canad. Math. Bull., {\bf 42}(1999), 285-290.

\bibitem{5} E. Fouvry, E. Kowalski and P. Michel, {\it Algebraic trace
functions over the primes}, available at
http://arxiv.org/abs/1211.6043v2, to appear in Duke Math. J.

\bibitem{6} E. Fouvry and P. Michel, {\it Sur certaines sommes
d'exponentielles sur les nombres premiers}, Ann. Sci. \'Ecole Norm.
Sup.(4), {\bf 31}(1998), no.1, 93-130.

\bibitem{7} E. Fouvry and I. E. Shparlinski, {\it On a ternary
quadratic form over primes}, Acta Arith., {\bf 150}(2011), 285-314.

\bibitem{8} D. Hajela, A. Pollington and B. Smith, {\it On Kloosterman
sums with oscillating coefficients}, Canad. Math. Bull., {\bf
31}(1988), 32-36.

\bibitem{9} G. Wang and Z. Zheng, {\it Kloosterman sums with oscillating
coefficients}, Chinese Ann. Math., {\bf 19}(1998), 237-242, in
Chinese; English translation in Chinese J. Contemp. Math., {\bf
19}(1998), 185-191.

\end{thebibliography}
\end{document}